\input amstex 
\documentstyle{amsppt}
\input bull-ppt
\keyedby{bull383/pxp}
\define\TVC{\operatorname{TVC}}

\topmatter
\cvol{28}
\cvolyear{1993}
\cmonth{April}
\cyear{1993}
\cvolno{2}
\cpgs{334-341}
\title Borel actions of Polish groups\endtitle
\author Howard Becker and Alexander S. Kechris\endauthor
\shortauthor{Howard Becker and A. S. Kechris}
\shorttitle{Borel actions of Polish groups}
\address Department of Mathematics, University of South 
Carolina,
Columbia, South Carolina 29208\endaddress
\ml becker\@cs.scarolina.edu\endml
\address Department of Mathematics, California Institute 
of Technology, 
Pasadena, California 91125\endaddress
\ml kechris\@romeo.caltech.edu\endml
\date April 16, 1992 and, in revised form, October 15, 
1992\enddate
\thanks The first author's research was partially 
supported by 
NSF Grant
DMS-8914426. The second author's research was partially 
supported by 
NSF Grant DMS-9020153\endthanks
\abstract We show that a {\it Borel\/} action of a Polish 
group on a standard Borel space
is Borel isomorphic to a {\it continuous\/} action of the 
group on a Polish
space, and we apply this result to three aspects of the 
theory of Borel actions
of Polish groups: universal actions, invariant probability 
measures, and the
Topological Vaught Conjecture. We establish the existence 
of universal actions
for any given Polish group, extending a result of Mackey 
and Varadarajan for the
locally compact case. We prove an analog of Tarski's 
theorem on paradoxical
decompositions by showing that the existence of an 
invariant Borel probability
measure is equivalent to the nonexistence of paradoxical 
decompositions with
countably many Borel pieces. We show that various natural 
versions of the
Topological Vaught Conjecture are equivalent with each 
other and, in the case
of the group of permutations of $\Bbb N$, with the 
model-theoretic Vaught
Conjecture for infinitary logic; this depends on our 
identification of the
universal action for that group.\endabstract
\endtopmatter

\document
A {\it Polish space\/} ({\it group\/}) is a separable, 
completely metrizable
topological space (group). A {\it standard Borel space\/} 
is a Polish space
with the associated Borel structure. A {\it Borel 
action\/} of a Polish group
$G$ on a standard Borel space $X$ is an action $(g,x)\in 
G\times X\mapsto
g_\cdot x$
of $G$ on $X$ which is Borel, as a function from the space 
$G\times X$ into
$X$. The structure of Borel actions of Polish locally 
compact, i.e., second
countable locally compact, topological groups has long 
been studied in ergodic
theory, operator algebras, and group representation 
theory. See, for example,
\cite{AM, Zi, Sin, VF, Mo, Ma1--Ma3, G, Var, FHM, Ra1, 
Ra2, K} for a sample of
works related to the themes that we will be studying here. 
More recently, there
has been increasing interest in the study of Borel actions 
of nonlocally compact
Polish groups. One instance is the {\it Vaught 
Conjecture\/}, a well-known open
problem in mathematical logic and the {\it Topological 
Vaught Conjecture\/}
(cf. \S\S1, 2) (see, e.g., \cite{Vau, Mi, St, Sa, L, Be, 
BM}). Another is the
ergodic theory and unitary group representation theory of 
so-called ``large
groups'' (see, e.g., \cite{Ve, O}). Also \cite{E} is 
relevant here.

Our purpose in this note is to announce a number of 
results about Borel actions
of general Polish groups. With the exception of Theorem 
2.1, these results are
new even for locally compact groups. Our Theorem 2.1 is 
known in the locally
compact case \cite{Ma2, Var}, but the proofs in this case 
relied on Haar
measure, so our proof seems new even in this case. The 
fundamental result is
Theorem 1.1, stating that for the actions considered Borel 
actions are
equivalent to continuous ones.

The paper is divided into three parts. The first part 
deals with Theorem 1.1
and a related result solving problems of Ramsay \cite{Ra2} 
 (raised in the
locally compact case) and Miller \cite{Mi} and, as a 
direct application,
shows the equivalence of three possible versions of the 
Topological Vaught
Conjecture. In the second part, we establish the existence 
of universal Borel
actions for any Polish group, extending a result of Mackey 
\cite{Ma2} and
Varadarajan \cite{Var} from the locally compact case. This 
is also applied to
establish the equivalence of the Topological Vaught 
Conjecture for the
symmetric group $S_\infty$, i.e., the permutation group of 
$\Bbb N$, with the
usual model-theoretic Vaught Conjecture for 
$L_{\omega_1\omega}$. The final
section deals with the problem of existence of invariant 
(countably additive)
Borel probability measures for a Borel action of a Polish 
group. It is a
well-known theorem of Tarski (see \cite{Wn}) that an 
arbitrary action of a group
$G$ on a set $X$ admits a {\it finitely additive\/} 
invariant probability
measure defined on all subsets of $X$ iff there is no 
``paradoxical
decomposition'' of $X$ with finitely many pieces. We show 
that there is a
complete analog of Tarski's theorem for Borel actions of 
Polish groups and
countably additive Borel invariant probability measures 
when we allow
``paradoxical decompositions'' to involve countably many 
Borel pieces. Our
proof uses the results of \S1 and the basic work of 
Nadkarni \cite{N}, who
proves this result in case $G=\Bbb Z$. 

\heading 1. Borel vs. topological group actions\endheading
By a {\it Borel\/} $G$-{\it space\/} we mean a triple 
$(X,G,\alpha)$, where $X$
is a standard Borel space, $G$ a Polish group, and 
$\alpha\:G\times X\to X$ a
Borel action of $G$ on $X$. We will usually write $X$ 
instead of $(X,G,\alpha)$
and $\alpha(g,x)=g_\cdot x$, when there is no danger of 
confusion. Two Borel
$G$-spaces $X$, $Y$ are {\it Borel isomorphic\/} if there 
is a Borel bijection
$\pi\:X\to Y$ with $\pi(g_\cdot x)=g_\cdot\pi(x)$. 
A {\it Polish\/} $G$-{\it space\/}
consists of a triple $(X,G,\alpha)$, where $X$ is a Polish 
space and
$\alpha\:G\times X\to X$ is continuous (thus any Polish 
$G$-space is also
Borel).

\thm{Theorem 1.1} Let $X$ be any Borel $G$-space. Then 
there is a Polish
$G$-space $Y$ Borel isomorphic to $X$.
\ethm

This answers a question of \cite{Ra2} (for locally compact 
$G)$ and \cite{Mi}.
Theorem 1.1 was known classically for discrete $G$ and was 
proved in \cite{Wh}
for $G=\Bbb R$. A convenient reformulation of our result 
is the following: If
$X$ is a Borel $G$-space, then there is a Polish topology 
$\tau$ on $X$ giving
the same Borel structure for which the action becomes 
continuous. In our proof,
sketched in \S4, we define $\tau$ explicitly and use a 
criterion of Choquet to
show that it is a Polish topology. One can also prove a 
version of Theorem 1.1
for more general ``definable'' actions of $G$ on separable 
metrizable spaces.

It is a classical result of descriptive set theory (see, 
e.g., \cite{Ku}) that
for any Polish space $X$ and any Borel set $B\subseteq X$, 
there is a Polish
topology, finer than the topology of $X$, and thus having 
the same Borel
structure as $X$, in which $B$ is clopen. We extend this 
result in the case of
Polish $G$-actions. The result below was known for 
$S_\infty$ (see \cite{Sa})
and is essentially a classical result for discrete $G$.

\thm{Theorem 1.2} Let $X$ be a Polish $G$-space and 
$B\subseteq X$ an invariant
Borel set. Then there is a Polish topology finer than the 
topology of $X$
\RM(and thus having the same Borel structure\RM) in which 
$B$ is now clopen and
the action is still continuous.
\ethm

The Topological Vaught Conjecture for a Polish group $G$, 
first conjectured by
Miller (see, e.g., \cite{Ro, p.~484}) has been usually 
formulated in one of
three, a priori distinct, forms. Let:

$\TVC_I(G)\Leftrightarrow$ For any Polish $G$-space $X$ 
and any
invariant Borel set $B\subseteq X$, either $B$ contains 
countably many orbits
or else there is a perfect set $P\subseteq B$ with any two 
distinct members of
$P$ belonging to different orbits;

$\TVC_{II}(G)\Leftrightarrow$ Same as $\TVC_I(G)$ but with 
$B=X$;

$\TVC_{III}(G)\Leftrightarrow$ For any Borel $G$-space $X$ 
either $X$
contains countably many orbits or else there is an 
uncountable Borel set
$P\subseteq X$ with any two distinct members of $P$ 
belonging to different
orbits.

Clearly, $\TVC_{III}(G)\Rightarrow \TVC_I(G)$, since a 
Borel set
in a Polish space is uncountable iff it contains a perfect 
set, and
$\TVC_I(G)\Rightarrow \TVC_{II}(G)$. From Theorem 1.1 we 
have 
$\TVC_{II}(G)\Rightarrow \TVC_{III}(G)$, thus

\thm{Corollary 1.3} For any Polish group $G$, all three 
forms of the
Topological Vaught Conjecture for $G$ are equivalent.
\ethm

The Topological Vaught Conjecture was motivated by the 
Vaught Conjecture in
logic, which we discuss in \S2. The truth or falsity of 
$\TVC(G)$ remains open
for general Polish $G$. It has been proved for $G$ locally 
compact (see
\cite{Sil}) or abelian (see \cite{Sa}). It can be also 
shown, using a method of
Mackey, that if $G$ is a closed subgroup of $H$, then 
$\TVC(H)$ implies
$\TVC(G)$. Of particular interest is the case $G=S_\infty$ 
(see \S2). It is
known (see \cite{Bu}) that for any Polish $G$-space $X$ 
there are either 
$\le\aleph_1$ many orbits or else perfectly many orbits, 
i.e., there is a
perfect set with any two distinct members of it belonging 
to different orbits.
So assuming the negation of the Continuum Hypothesis (CH), 
there are
$2^{\aleph_0}$ many orbits iff there are perfectly many 
orbits. Thus the
meaning of $\TVC$ is that the set of orbits cannot be a 
counterexample to CH.

\heading 2. Universal actions\endheading
Let $X$, $Y$ be Borel $G$-spaces. A {\it Borel 
embedding\/} of $X$ into $Y$ is
a Borel injection $\pi\:X\to Y$ such that $\pi(g_\cdot 
x)=g_\cdot 
\pi(x)$. Note that
$\pi[X]$ is an invariant Borel subset of $Y$. By the usual 
Schroeder-Bernstein
argument, $X$, $Y$ can be Borel-embedded in each other iff 
they are Borel
isomorphic. A Borel $G$-space $\scr U$  is {\it 
universal\/} if every Borel
$G$-space $X$ can be Borel-embedded into $\scr U$. It is 
unique up to Borel
isomorphism.

\thm{Theorem 2.1} For any Polish group $G$, there is a 
universal Borel
$G$-space $\scr U_G$. Moreover, $\scr U_G$ can be taken to 
be a Polish
$G$-space.
\ethm

Actually the proof of Theorem 2.1 shows that one can 
Borel-embed in $\scr U_G$
any Borel action of $G$ on a separable metrizable space. 
In particular, a Borel
action of $G$ on an analytic Borel space (i.e., an 
analytic set with its
associated Borel structure) is Borel isomorphic to a 
continuous action of $G$
on an analytic space, i.e., an analytic set with its 
associated topology. For
locally compact $G$, Theorem 2.1 has been proved in 
\cite{Ma2, Var}. In this
case, $\scr U_G$ can be taken to be compact. Our proof 
gives a new proof of
this result, with a different universal space, which 
avoids using the Haar
measure. It is unknown whether $\scr U_G$ can be taken to 
be compact in the
general case.

In the particular case of the group $G=S_\infty$, with the 
Polish topology it
inherits as a $G_\delta$ subspace of the Baire space $\Bbb 
N^{\Bbb N}$, our
proof of Theorem 2.1 gives a particularly simple form of 
$\scr U_G$, which we
use below. Let $X_\infty$ be the space of all maps from 
the set $\Bbb N^{<
\Bbb N}$ of finite sequences of natural numbers into 
$2=\{0,1\}$. Clearly this
space is homeomorphic to the Cantor space. Consider the 
following action of
$S_\infty$ on $X_\infty$:
$g_\cdot x(s_0,s_1,\dotsc,s_{n-1})=x(g^{-1}(s_0),%
\dotsc,g^{-1}(s_{n-1}))$, if
$s=(s_0,s_1,\dotsc,s_{n-1})\in\Bbb N^{<\Bbb N}$. Then 
$X_\infty$ with this
action is a universal $S_\infty$-space and is clearly 
compact. We would like
to thank R. Dougherty for suggesting that we can use the 
above space instead of
a somewhat more complicated one we had originally.

This leads to an interesting fact concerning the Vaught 
Conjecture. The Vaught
Conjecture for $L_{\omega_1\omega}$ is the assertion that 
every $L_{\omega_1
\omega}$ sentence in a countable language has either 
countably many or else
perfectly many countable models, up to isomorphism. Let 
$L$ be a countable
language which we assume to be relational, say 
$L=\{R_i\}_{i\in I}$, where $I$
is a countable set and $R_i$ is a $n_i$-ary relation 
symbol. Let
$X_L=\prod_{i\in I}2^{\Bbb N^{n_i}}$ which is homeomorphic 
to the Cantor space
(if $L\ne \emptyset)$. 
We view $X_L$ as the space of countably infinite 
structures for $L$, 
identifying $x=(x_i)_{i\in I}\in X_L$ with the structure
$\scr A=\langle\Bbb N$, $R_i^{\scr A}\rangle$, where 
$R_i^{\scr
A}(s)\Leftrightarrow x_i(s)=1$. The group $S_\infty $ acts 
in the obvious way
on $X_L\:g_\cdot x=y\Leftrightarrow\forall 
i[y_i(s_0,\dotsc,s_{n_i-1})=1$ iff 
$x_i(g^{-1}(s_0),\dotsc,g^{-1}(s_{n_i-1}))=1]$. Thus if 
$x$, $y$ are identified 
with the structures $\scr A$, $\scr B$ resp., $g_\cdot 
x=y$ iff $g$ is an isomorphism
of $\scr A$, $\scr B$. This action is called the {\it 
logic action\/} of
$S_\infty$ on $X_L$. By \cite{L-E} the Borel  invariant 
subsets of this action
are exactly the sets of models of $L_{\omega_1\omega}$ 
sentences. So the
Vaught Conjecture for $L_{\omega_1\omega}$ (VC), a 
notorious open problem, is
the assertion that for any countable $L$ and any Borel 
invariant subset
$B\subseteq X_L$, either $B$ contains countably many 
orbits or else there is a
perfect set $P\subseteq B$, no two distinct members of 
which are in the same
orbit. It is thus a special case of the Topological Vaught 
Conjecture for
$S_\infty$. (Historically, of course, $\TVC$ came much 
later and was inspired
by VC.) The universal Borel $S_\infty$-space $X_\infty$ 
above is the same as
$X_{L_\infty}$ for $L_\infty=\{R_n\}_{n\in\Bbb N}$, with 
$R_n$ an $n$-ary
relation, so we have

\thm{Corollary 2.2} The logic action on the space of 
structures of the language
containing an $n$-ary relation symbol for each $n\in\Bbb 
N$ is a universal
$S_\infty$-space. In particular, Vaught{\rm '}s conjecture 
for $L_{\omega_1
\omega}$ is equivalent to the Topological Vaught 
Conjecture for $S_\infty$.
\ethm

We conclude with an application to equivalence relations. 
Given two equivalence
relations $E$, $F$ on standard Borel spaces $X$, $Y$ 
resp., we say that $E$ is
{\it Borel embeddable\/} in $F$ iff there is a Borel 
injection $f\:X\to Y$ with
$xEy\Leftrightarrow f(x)Ff(y)$. Given a class of 
equivalence relations $\scr S$,
a member $F\in\scr S$ is called {\it universal\/} if every 
$E\in\scr S$ is
Borel embeddable in $F$. For each Borel $G$-space denote 
by $E_G$ the
corresponding equivalence relation induced by the orbits 
of the action:
$xE_Gy\Leftrightarrow \exists g\in G(g_\cdot x=y)$. The 
following is based on the
proof of Theorem 2.1 and the work in \cite{U}, a paper 
which was brought to our
attention by W. Comfort.

\thm{Corollary 2.3} There is a universal equivalence 
relation in the class of
equivalence relations $E_G$ induced by Borel actions of 
Polish groups.
\ethm

\heading 3. Invariant measures\endheading
Let $G$ be a group acting on a set $X$. Given 
$A,B\subseteq X$, we say that
$A$, $B$ are {\it equivalent by finite decomposition\/}, 
in symbols 
$A\sim B$, if
there are partitions $A=\bigcup_{i=1}^nA_i$, 
$B=\bigcup_{i=1}^nB_i$, and
$g_1,\dotsc,g_n$ with $g_i.A_i=B_i$. We say that $X$ is 
$G$-{\it paradoxical\/}
if $X\sim A\sim B$ with $A\cap B=\emptyset$. A {\it 
finitely additive
probability\/} (f.a.p.) measure on $X$ is a map 
$\varphi\:\scr P(X)\to[0,1]$
such that $\varphi(X)=1$, $\varphi(A)+
\varphi(B)=\varphi(A\cup B)$, if
$A\cap B=\emptyset$.  Such a  $\varphi$ is $G$-{\it 
invariant\/} if $\varphi(A)=
\varphi(g_\cdot A)$ for all $g\in G$ and $A\subseteq X$. 
If $X$ is $G$-paradoxical,
there can be no $G$-invariant f.a.p. measure on $X$. A 
well-known theorem of
Tarski asserts that the converse is also true (see 
\cite{Wn}).

\thm{Theorem 3.1 {\rm (Tarski)}\rm} Let a group $G$ act on 
a set $X$. Then
there is a $G$-invariant finitely additive probability 
measure on $X$ iff $X$
is not $G$-paradoxical.
\ethm

It is natural to consider to what extent Tarski's theorem 
goes through for
countably additive probability measures. Let $(X,\scr A)$ 
be a measurable
space, i.e., a set equipped with a $\sigma$-algebra. Let a 
group $G$ act on $X$
so that $A\in\scr A\Rightarrow g_\cdot 
A\in\scr A$. Given $A,B\in\scr A$, we say that
$A$, $B$ are {\it equivalent by countable 
decomposition\/}, in symbols
$A\sim_\infty B$, if there are partitions $A=\bigcup_{i\in 
I}A_i$,
$B=\bigcup_{i\in I}B_i$, with $I$ countable and 
$A_i,B_i\in\scr A$, and
$\{g_i\}_{i\in I}$ so that $g_{i^\cdot}A_i=B_i$. We say 
that $(X,\scr
A)$, or just $X$ if there is no danger of confusion, is 
{\it countably\/}
$G$-{\it paradoxical\/} if $X\sim_\infty A\sim_\infty B$ 
with $A,B\in\scr A$
and $A\cap B=\emptyset$. A probability measure $\mu$ on 
$(X,\scr A)$ is
$G$-{\it invariant\/} if $\mu(A)=\mu(g_\cdot A)$ for all 
$g\in G$ and $A\in\scr A$.
Again, if there is a $G$-invariant probability measure on 
$(X,\scr A)$, $X$
cannot be countably $G$-paradoxical. Is the converse true? 
It turns out that
the answer is negative in this generality. See \cite{Wn, 
Za} for more on the
history and some recent developments on this problem. We 
show here, however,
that the converse holds in most regular situations, i.e., 
when $G$ is Polish
and acts in a Borel way on a standard Borel space $X$.

\thm{Theorem 3.2} Let $X$ be a Borel $G$-space. Then the 
following are
equivalent\RM:
\roster
\item \<$X$ is not countably $G$-paradoxical\RM;
\item there is a $G$-invariant Borel probability measure 
on $X$.
\endroster
\ethm

The proof of this theorem is based on the results in \S1 
and \cite{N}, which
proves Theorem 3.2 for the case $G=\Bbb Z$. It turns out 
that Theorem 3.2 holds
also for any {\it continuous\/} action of a separable 
topological group on a
Polish space $X$.

\heading 4.  Sketches of proofs\endheading
\demo{For Theorem {\rm 1.1}} We have a Borel action of $G$ 
on $X$. Fix a
countable basis $\scr B$ for the topology of $G$. For 
$A\subseteq X$ and
$U\subseteq G$ open, let $A^{\Delta U}=\{x:g_\cdot x\in A$ 
for a set of $g$'s
which is nonmeager in $U\}$ be the {\it Vaught 
transform\/} \cite{Vau}. We first
find a countable Boolean algebra $\scr C$ of Borel subsets 
of $X$ such that
(1) $A\in \scr C\Rightarrow A^{\Delta U}\in \scr C$ for 
$U\in\scr B$, and (2)
the topology generated by $\scr C$ is Polish. Then let 
$\tau$ be the topology
on $X$ generated by $\{A^{\Delta U}\:A\in \scr C$, 
$U\in\scr B\}$. It suffices
to show that $\tau$ is a Polish topology and that the 
action of $G$ is
continuous in this topology. Then since $\tau$ consists of 
Borel sets in $X$
and gives a Polish topology, it gives rise to the original 
Borel structure of
$X$.

The continuity of the action may be checked by direct 
computation. To see that
$\tau$ gives a Polish topology, we first check that it is 
$T_1$ and regular, as
well as obviously second countable, hence metrizable, and 
then apply a
criterion of Choquet to conclude that it is Polish. 
Choquet's criterion may be
stated as follows.

Associate to $X$ the {\it strong Choquet game\/}, in which 
the first player
specifies a sequence of open sets $U_n$ and elements 
$x_n\in U_n$, while the
second player responds with open sets $V_n$, satisfying 
$x_n\in V_n\subseteq
U_n$ and $U_{n+1}\subseteq V_n$; the second player wins if 
the intersection
of the $U_n$ is nonempty. The space $X$ is called {\it a 
strong Choquet
space\/} (see, e.g., \cite{HKL}) 
if the second player has a winning strategy for this game, 
and
Choquet's criterion (unpublished, but see \cite{C} for a 
related version)
states that a topological space is Polish if and only if 
it is separable,
metrizable, and strong Choquet.
\enddemo

\demo{For Theorem {\rm 1.2}} The proof is quite similar, 
taking a somewhat
larger $\scr C$ and including a basis for the original 
topology in $\tau$.
\enddemo

\demo{For Theorem {\rm 2.1}} Let $\scr F(G)$ be the space 
of closed subsets of
$G$ with the (standard) Effros Borel structure, i.e., the 
one generated by the
sets $\scr F_V=\{F\in \scr F(G)\:F\cap V\ne \emptyset\}$, 
for $V\subseteq G$ open.
Let $G$ act on $\scr F(G)$ by left multiplication and let 
$\scr U_G=\scr F(G)^{
\Bbb N}$, with $G$ acting coordinatewise. Given any Borel 
$G$-space $X$, let
$\{S_n\}$ be a sequence of Borel subsets of $X$ separating 
points. For each 
$A\subseteq G$, let $D(A)=\{g\in G$: For every 
neighborhood $U$ of $g$, $A$ is
nonmeager in $U\}$ and map $x\in X$ into the sequence 
$f(x)=\{D(\widetilde
S_n)^{-1}\}\in \scr U_G$, where $\widetilde 
S_n=\{g\:g_\cdot x\in S_n\}$. This
is an embedding of $X$ into $\scr U_G$. We can now make 
$\scr U_G$ into a
Polish $G$-space using Theorem 1.1. When $G$ is locally 
compact and $\overline
G=G\cup\{\infty\}$ is the one-point compactification of 
$G$, then we can extend
the action of $G$ on itself by left multiplication to 
$\overline G$ by setting
$g_\cdot\infty=\infty$. Then instead of $\scr F(G)$ we can 
use $\scr
F(\overline G)=K(\overline G)=$the compact metrizable 
space of compact subsets
of $\overline G$, with the obvious $G$-action, and the 
universal space is now
$K(\overline G)^{\Bbb N}$, which is compact. 
(Alternatively, we can use the
Fell topology on $\scr F(G)$.)
\enddemo

\demo{For Corollary {\rm 2.2}} The action of $S_\infty$ on 
itself by left
multiplication extends to the action of $S_\infty$ on 
$\Bbb N^{\Bbb N}$ (the
Baire space) by left composition. So instead of $\scr 
F(S_\infty)$, one can use
$\scr F(\Bbb N^{\Bbb N})$ and thus can take $\scr F(\Bbb 
N^{\Bbb N})^{\Bbb N}$
as a universal space. But closed subsets of $\Bbb N^{\Bbb 
N}$ can be identified
with trees on $\Bbb N$, i.e., subsets of $\Bbb N^{<\Bbb 
N}$ closed under initial
segments. A tree can be viewed now as a sequence of 
$n$-ary relations on $
\Bbb N$, the $n$\<th relation identifying which $n$-tuples 
belong to the tree.
Thus $\scr F(\Bbb N^{\Bbb N})^{\Bbb N}$ can be embedded in 
the space 
of structures
$X_L$, for the language containing infinitely many $n$-ary 
relation symbols for
each $n$, and then by a simple coding pointed out by R. 
Dougherty, in the space
of structures $X_{L_\infty}$ for the language $L_\infty$.
\enddemo

\demo{For Theorem {\rm 3.2}} One first shows that the 
result holds for countable
groups by adapting appropriately the proof of \cite{N}. 
For an arbitrary $G$,
let $G_1$ be a countable dense subgroup of $G$ and notice 
that if $G$ acts
continuously, then any $G_1$-invariant measure is also 
$G$-invariant. Then
apply Theorem 1.1.
\enddemo

We would like to thank A. Louveau for many helpful 
suggestions concerning the
results announced in this paper and G. Cherlin, who helped 
in improving
considerably the exposition.

\enddocument